\documentclass[11pt,reqno]{amsart}

\usepackage{amsmath,amsthm}
\usepackage{amssymb}
\usepackage{euscript}

\numberwithin{equation}{subsection}

\newtheorem{Thm}[subsection]{Theorem}

\newtheorem{Prop}[subsection]{Proposition}

\theoremstyle{definition}

\begin{document}
\title[]{Vector bundles over iterated suspensions of stunted real projective spaces}
\author{Aniruddha C. Naolekar}\address{Indian Statistical Institute, 8th Mile, Mysore Road, RVCE Post, Bangalore 560059, INDIA.}
\email{ani@isibang.ac.in}
\author{ Ajay Singh Thakur}\address{School of Maths, TIFR, Homi Bhabha Road, Colaba, Mumbai 400005, INDIA.}
\email{athakur@math.tifr.res.in}

\begin{abstract}
Let $X^k_{m,n}=\Sigma^k (\mathbb R\mathbb P^m/\mathbb R\mathbb P^n)$. In this note we completely determine the values of $k,m,n$ for which the total Stiefel-Whitney class $w(\xi)=1$ for any vector bundle $\xi$ over $X^k_{m,n}$.  
\end{abstract}


\email{}

\date{}
\maketitle

\section{Introduction}

Recall (see \cite{tanaka1}) that a $CW$-complex $X$ is said to be $W$-trivial if for every vector bundle $\xi$ over $X$ the total Stiefel-Whitney class $w(\xi)=1$. A theorem of Atiyah-Hirzebruch (\cite{atiyah}, Theorem 2) says that for any finite $CW$-complex $X$, the $9$-fold suspension $\Sigma^9X$ is $W$-trivial. In the same paper Atiyah-Hirzebruch have shown (\cite{atiyah}, Theorem 1) that the sphere $S^k$ is $W$-trivial if and only if $k\neq 1,2,4,8$ (see also, \cite{milnor}, Theorem 1).  

In view of the Atiyah-Hirzebruch theorem  it is interesting to understand whether or not the iterated suspension $\Sigma^kX$ of a finite $CW$-complex $X$ is $W$-trivial with $0\leq k\leq 8$. In recent times there has been some interest in understanding $W$-triviality of iterated suspensions of spaces (see, \cite{tanaka1}, \cite{tanaka3}, \cite{ajay} and the references therein).  

In \cite{tanaka3}, the author has completely determined the values of $k$ and $n$ for which the iterated suspension $\Sigma^k \mathbb F\mathbb P^n$ is $W$-trivial. Here $\mathbb F\mathbb P^n$ denotes the projective space of $1$-dimensional subspaces of $\mathbb F^{n+1}$ where $\mathbb F$ is the field $\mathbb R$ of reals, the field $\mathbb C$ of complex numbers or the skew-field $\mathbb H$ of quaternions. In \cite{tanaka1}, the author has completely described the cases under which the stunted projective space $\mathbb R\mathbb P^m/\mathbb R\mathbb P^n$ is $W$-trivial. In \cite{ajay}, the second author has almost complete results concerning $W$-triviality of the iterated suspension $\Sigma^kD(m,n)$ of the Dold manifold $D(m,n)$. 

Let $X_{m,n}$ denote the stunted projective space $\mathbb R\mathbb P^m/\mathbb R\mathbb P^n$ and let $X^k_{m,n}$ denote the $k$-fold suspension 
$\Sigma^kX_{m,n}$ of $X_{m,n}$. In this note we completely determine the values of $k,m,n$ for which $X^k_{m,n}$ is $W$-trivial. 

In view of the Atiyah-Hirzebruch theorem, we assume that $0\leq k\leq 8$. Also note that the cases $X^k_{m,0}=\Sigma^k\mathbb R\mathbb P^m$ and 
$X^k_{m,m-1}=\Sigma^kS^m=S^{m+k}$ are completely understood. In the sequel we assume that $0<n<m$ and hence, in particular, $m\geq 2$. Since the case $X_{m,n}$ is completely understood, we state our results for $X^k_{m,n}$ with $1\leq k\leq 8$. 
The following statements completely describe the cases when $X^k_{m,n}$ is $W$-trivial.  

\begin{Thm} \label{maintheorem} Let $X^k_{m,n}$ be as above with $0<n<m$. \begin{enumerate}
\item If $k=3,5$, then $X^k_{m,n}$ is $W$-trivial if and only if $m+k\neq 8$. 
\item $X^6_{m,n}$ is $W$-trivial if and only if $(m,n)\neq (3,1), (2,1)$. 
\item $X^7_{m,n}$ is $W$-trivial. 
\end{enumerate}

\end{Thm}

\begin{Thm}\label{secondtheorem} 
Let $0<n<m$. Then $X^1_{m,n}$ is $W$-trivial if and only if $m\neq 3,7$. 
\end{Thm}


\begin{Thm} \label{thirdtheorem} Let $0<n<m$. Then 
\begin{enumerate}
\item $X^2_{m,n}$ is not $W$-trivial if $n=1$. 
\item $X^2_{m,n}$ is $W$-trivial if $m\not\equiv 0,6,7 \mbox{ {\rm (mod $8$)}}$ and $n\geq 2$. 
\item $X^2_{8t+6,n}$ is $W$-trivial if $t\geq 1$ and $n\geq 2$. $X^2_{6,n}$ is not $W$-trivial. 
\item $X^2_{8t+7,n}$ is $W$-trivial if $t\geq 1$ and $n\geq 2$. $X^2_{7,n}$ is $W$-trivial if and only if $n=6$.
\item $X^2_{8t,n}$ is $W$-trivial if $n\geq 2$.
\end{enumerate} 
\end{Thm}

\begin{Thm}\label{fourththeorem}
Let $0<n<m$. Then  
\begin{enumerate}
\item $X^4_{m,n}$ is $W$-trivial if $m=2,3$ and not $W$-trivial if $m=4$, 
\item $X^4_{m,n}$ is not $W$-trivial if $m>4$ and $n=1,2,3$. $X^4_{m,n}$ is $W$-trivial if $m>4$ and $n\geq 4$.
\end{enumerate}  
\end{Thm}

\begin{Thm}\label{fifththeorem}
Let $0<n<m$. Then $X^8_{m,n}$ is $W$-trivial. 
\end{Thm}

The proofs of the above theorems crucially make use of the computations of the $\widetilde{KO}$-groups of the real projective spaces and the 
stunted real projective spaces. 
In the next section we first state some easy to verify general observations and prove our main results. 

{\em Conventions.} All references to cohomology groups will mean singular cohomology with $\mathbb Z_2$-coefficients. Given a map $\alpha:X\longrightarrow Y$, the induced homomorphism in cohomology and $\widetilde{KO}$-groups will again be denoted by $\alpha$. 


\section{Proof of the theorems}

We begin by recording some important observations which will be crucial in the proofs of the main theorems. 

\begin{Prop} \label{firstlemma}\begin{enumerate}
\item $X^k_{m,n}$ is $W$-trivial if there does not exist an integer $s$ such that $n+1+k\leq 2^s\leq m+k$. 
\item Suppose $0\leq n'<n<m$. If $X^k_{m,n'}$ is $W$-trivial, then so is $X^k_{m,n}$. In particular, if $\Sigma^k \mathbb R\mathbb P^m$ is  $W$-trivial, then so is $X^k_{m,n}$. 
\end{enumerate}\end{Prop}
{\em Proof.} (1) follows from the well-known fact that the first non-zero Stiefel-Whitney class of a vector bundle is in degree a power of $2$. 
(2) follows from the fact that the obvious map $ X^k_{m,n'}\longrightarrow X^k_{m,n}$ induces isomorphism in cohomology in degree $i$ with $n+1+k\leq i\leq m+k$. \qed

We note that if $m$ is odd, then we have a splitting 
$$X_{m, m-2}=S^m\vee S^{m-1}$$
and if $m$ is even, then we have 
$$X_{m,m-2}= M(\mathbb Z_2, m-1)=\Sigma^{m-2}\mathbb R\mathbb P^2$$
where $M(\mathbb Z_2, m-1)$ denotes the Moore space of type $(\mathbb Z_2, m-1)$. Also note that if $\widetilde{KO}^{-k}(X)=0$, then $\Sigma^kX$ is $W$-trivial. 

Given a sequence of integers $0\leq p<n<m$, the cofiber sequence 
$$X_{n,p}\longrightarrow X_{m,p}\longrightarrow X_{m,n}$$
gives rise to an exact sequence 
$$\cdots\rightarrow\widetilde{KO}^{-i}(X_{m,n})\longrightarrow \widetilde{KO}^{-i}(X_{m,p})\longrightarrow\widetilde{KO}^{-i}(X_{n,p})\longrightarrow\widetilde{KO}^{-i+1}(X_{m,n})\rightarrow\cdots.$$
Our proofs, in many cases, involve analyzing the above exact sequence of $\widetilde{KO}$-groups corresponding to a suitable choice of a cofiber sequence as above. 
We shall use the above observations implicitly in the sequel.

We record here what is known about the $W$-triviality of $\Sigma^k\mathbb R\mathbb P^m$ and the stunted projective spaces $X_{m,n}$ for easy reference.

\begin{Thm}\label{tanakastheorem}{\rm (\cite{tanaka3}, Theorem 1.4)} (1) If $k=1,2,4,8$, then $\Sigma^k\mathbb R\mathbb P^m$ is not $W$-trivial if and only if $m\geq k$. (2) If $k=3,5,7$, then 
$\Sigma^k\mathbb R\mathbb P^m$ is not $W$-trivial if and only if $m+k=4,8$. (3) $\Sigma^6\mathbb R\mathbb P^m$ is not $W$-trivial if and only if $m=2,3$. \qed
\end{Thm}

\begin{Thm}\label{tanakastheorem1}{\rm(\cite{tanaka1}, Theorem 4.1)} 
Suppose $1\leq n \leq m-2$. Then $X_{m,n}$ is $W$-trivial if and only if $m<2^{\varphi(n)}$ where $\varphi(n)$ denotes the number of integers $i$ 
such that $0<i\leq n$ and $i\equiv 0,1,2,4 \mbox{ {\rm (mod $8$)}}$.\qed
\end{Thm}

{\em Proof of Theorem\,\ref{maintheorem}}. 
Assume first that $k=3,5$. If $m+k\neq 8$ then as $\Sigma^k\mathbb R\mathbb P^m$ is $W$-trivial it follows 
from Proposition\,\ref{firstlemma} (2) that $X^k_{m,n}$ is also $W$-trivial. 
Next we look at $X^3_{5,n}$. The obvious map $X^3_{5,n}\longrightarrow S^8$ induces isomorphism in the top cohomology and hence the Hopf bundle on $S^8$ pulls back to a bundle $\xi$ over $X^3_{5,n}$ with $w(\xi)\neq 1$. A similar argument works for $X^5_{3,n}$. This completes the proof of (1).


The case (2) when $m\neq 2,3$ follows from arguments similar to the above case. That $X^6_{3,n}$ is not $W$-trivial for $n=1,2$ follows as from 
the facts that 
$X^6_{3,1}=S^9\vee S^8$ and $X^6_{3,2}=S^8$. Clearly, $X^6_{2,1}=S^8$ is not $W$-trivial. This completes the proof of (2). 

Finally as $m\geq 2$, $W$-triviality of $\Sigma^7\mathbb R\mathbb P^m$ implies that $X^7_{m,n}$ is $W$-trivial. This completes the proof of (3) and the theorem. 


{\em Proof of Theorem\,\ref{secondtheorem}}. 
There are obvious maps $X^1_{3,n}\longrightarrow S^4$ and $X^1_{7,n}\longrightarrow S^8$ that induce isomorphisms in cohomology in the top dimension. The Hopf bundles on $S^4$ and $S^8$ then pull back to give bundles with total Stiefel-Whitney class not equal to $1$. This shows that $X^1_{m,n}$ is not $W$-trivial if $m=3,7$. 

Next assume that $m\not\equiv 3,7 \mbox{ {\rm (mod $8$)}}$. Then by (\cite{fujii}, Table 3) we have $\widetilde{KO}^{-1}(X_{m,1})=0$. Hence when 
$m\not\equiv 3,7 \mbox{ {\rm (mod $8$)}}$, we have that $X^1_{m,1}$ and hence $X^1_{m,n}$ is $W$-trivial for all $n\geq 1$. 

Finally we look at the case when $m=8t+3, 8t+7$ with $t\geq 1$. 
First consider $X^1_{m,n}$ with $m=8t+3$ and $t\geq 1$. Consider the exact sequence 
$$\cdots\rightarrow\widetilde{KO}^{-1}(X_{8t+3,8t+2})\stackrel{\alpha}\longrightarrow\widetilde{KO}^{-1}(X_{8t+3,1})\rightarrow\widetilde{KO}^{-1}
(X_{8t+2,1})\rightarrow\cdots.$$
As $8t+2\not\equiv 3,7 \mbox{ {\rm (mod $8$)}}$, the last group in the above sequence is zero (\cite{fujii}, Table 3) and hence $\alpha$ is an epimorphism. 
Hence $W$-triviality of the sphere $X^1_{8t+3,8t+2}$ implies $W$-triviality of $X^1_{8t+3,1}$. Hence $X^1_{m,n}$ is $W$-trivial for all $n\geq 1$ when 
$m=8t+3$ and $t\geq 1$. The case $X^1_{8t+7,1}$ is dealt with similarly by looking at the cofiber sequence 
$$X_{8t+6,1}\longrightarrow X_{8t+7,1}\longrightarrow X_{8t+7,8t+6}.$$
Thus $X^1_{m,n}$ is $W$-trivial if and only if $m\neq 3,7$. This completes the proof of the theorem.

{\em Proof of Theorem\,\ref{thirdtheorem}}. We prove each of the claims in the theorem. 

{\em Proof of (1)}. We consider the exact sequence 
$$\cdots \rightarrow\widetilde{KO}^{-2}(X_{m,2})\longrightarrow \widetilde{KO}^{-2}(X_{m,1})\stackrel{\alpha}\longrightarrow\widetilde{KO}^{-2}(S^2)\longrightarrow 
\widetilde{KO}^{-1}(X_{m,2})\rightarrow\cdots. $$
 Assume that $m\not\equiv 3,7 \mbox{ {\rm (mod $8$)}}$. It then follows from (\cite{fujii}, Table 4) that the last group in the above sequence is zero. Hence the homomorphism $\alpha$ is an epimorphism. Thus there exists a vector bundle $\xi$ over $X^2_{m,1}$ that pulls back to the Hopf bundle over $S^4$. Hence $w(\xi)\neq 1$ showing that if $m\not\equiv 3,7 \mbox{ {\rm (mod $8$)}}$, then $X^2_{m,1}$ is not $W$-trivial. 
 
 We now deal with the cases $m\equiv 3,7 \mbox{ {\rm (mod $8$)}}$. 
Consider the exact sequence 
$$\cdots \rightarrow\widetilde{KO}^{-2}(X_{8t+3,8t+2})\rightarrow \widetilde{KO}^{-2}(X_{8t+3,1})\stackrel{\alpha}\rightarrow\widetilde{KO}^{-2}(X_{8t+2,1})\rightarrow 
\widetilde{KO}^{-1}(X_{8t+3,8t+2})\stackrel{\beta}\rightarrow$$
$$\rightarrow\widetilde{KO}^{-1}(X_{8t+3,1})\rightarrow\widetilde{KO}^{-1}(X_{8t+2,1})\rightarrow\cdots. $$
The last group in the above sequence is zero (\cite{fujii}, Table 3). Now $\widetilde{KO}^{-1}(X_{8t+3,8t+2})=\mathbb Z$ and 
$ \widetilde{KO}^{-1}(X_{8t+3,1}) = \mathbb Z$ by (\cite{fujii}, Table 3). Thus $\beta$ is an isomorphism and hence $\alpha$ is an epimorphism. As 
$8t+2\not\equiv 3, 7 \mbox{ {\rm (mod $8$)}}$ it follows that $X^2_{8t+2,1}$ is not $W$-trivial and hence $X^2_{8t+3,1}$ is not $W$-trivial. The case $X^2_{8t+7,1}$ is dealt with similarly by looking at the cofiber sequence  $X_{8t+6,1}\rightarrow X_{8t+7,1}\rightarrow X_{8t+7,8t+6}$. This completes the proof of (1). 

{\em Proof of (2)}. Assume that $m\not\equiv 0,6,7 \mbox{ {\rm (mod $8$)}}$. By (\cite{fujii}, Table 4), $\widetilde{KO}^{-2}(X_{m,2})=0$. Thus $X^2_{m,n}$ is $W$-trivial for all $n\geq 2$ if $m\not\equiv 0,6,7 \mbox{ {\rm (mod $8$)}}$. This completes the proof of (2). 

{\em Proof of (3)}. We consider $X^2_{8t+6,n}$. If $t=0$, the obvious map $X^2_{6,n}\longrightarrow S^8$ induces isomorphism in cohomology in the top dimension and hence 
the Hopf bundle on $S^8$ pulls back to a bundle $\xi$ with $w(\xi)\neq 1$. 

Assume now that $t\geq 1$. In the exact sequence 
$$\cdots \rightarrow\widetilde{KO}^{-2}(X_{8t+6,8t+5})\stackrel{\alpha}\longrightarrow \widetilde{KO}^{-2}(X_{8t+6,2})\longrightarrow\widetilde{KO}^{-2}(X_{8t+5,2})\rightarrow \cdots$$
$\alpha$ is an epimorphism as the last group is zero (\cite{fujii}, Table 4). Clearly $X^2_{8t+6,8t+5}$ is $W$-trivial showing that $X^2_{8t+6,2}$ is $W$-trivial. Hence $X^2_{8t+6,n}$ is $W$-trivial for all $n\geq 2$. This completes the proof of (3).  

{\em Proof of (4)}. First note that as $X^2_{7,5}=S^9\vee S^8$ we have that $X^2_{7,5}$ is not  $W$-trivial. Thus  $X^2_{7,n}$ is not $W$-trivial for $n\leq 5$. Clearly, $X^2_{7,6}=S^9$ is $W$-trivial.  

Next we consider $X^2_{8t+7,n}$ with $t\geq 1$. In the exact sequence 
$$\rightarrow\widetilde{KO}^{-2}(X_{8t+7,8t+6})\rightarrow \widetilde{KO}^{-2}(X_{8t+7,2})\stackrel{\alpha}\rightarrow\widetilde{KO}^{-2}(X_{8t+6,2})\rightarrow 
\widetilde{KO}^{-1}(X_{8t+7,8t+6})\rightarrow$$
 as $ \widetilde{KO}^{-2}(X_{8t+6,2})=\mathbb Z_2$ (\cite{fujii}, Table 4) and the last group is infinite cyclic we have that $\alpha$ is an epimorphism. By (3) above, $X^2_{8t+6,2}$ is $W$-trivial and hence $X^2_{8t+7,2}$ is $W$-trivial. Hence the proof of (4) is complete. 

{\em Proof of (5)}. We first concentrate on the case $t=1$. Note that $X^2_{8,7}=S^{10}$ and hence is $W$-trivial. As 
$$X^2_{8,6}=\Sigma^2M(\mathbb Z_2,7)=\Sigma^8\mathbb RP^2,$$
we have, by Theorem\,\ref{tanakastheorem}, that $X^2_{8,6}$ is $W$-trivial. Next consider the exact sequence 
$$\rightarrow\widetilde{KO}^{-2}(X_{8,6})\rightarrow\widetilde{KO}^{-2}(X_{8,5})\stackrel{\beta}\longrightarrow\widetilde{KO}^{-2}(S^6)\stackrel{\alpha}\rightarrow\widetilde{KO}^{-1}(X_{8,6})\rightarrow\widetilde{KO}^{-1}(X_{8,1})\rightarrow.$$
We note that $\widetilde{KO}^{-2}(S^6)=\mathbb Z$ generated by the Hopf bundle $\nu$, $\widetilde{KO}^{-1}(X_{8,6})=\mathbb Z_2$ (\cite{fujii}, Table 4),   $\widetilde{KO}^{-1}(X_{8,1})=0$ (\cite{fujii}, Table 3) and $\widetilde{KO}^{-2}(X_{8,5})=\mathbb Z\oplus\mathbb Z_2$ (\cite{fujii}, Table 3). Hence $\alpha $ is an epimorphism. 
Since $\beta$ induces isomorphism in $8th$ cohomology, the equality $w(\beta(\xi))=1$ implies that $w(\xi)=1$. 
It follows from the exactness of the above sequence that there is a generator $\xi$ of the torsion free part of $\widetilde{KO}^{-2}(X_{8,5})$ with  $\beta(\xi)=2\nu$. Since $w(2\nu)=1$, we have that $w(\xi)=1$. If $\eta$ is the generator of the torsion part, $w(\beta(\eta))=1$ as $\beta(\eta)$ is (stably) trivial. Hence $w(\eta)=1$. 
This shows that for any $\theta\in \widetilde{KO}^{-2}(X_{8,5})$, $w(\theta)=1$ 
and hence $X^2_{8,5}$ is $W$-trivial. 
$W$-triviality of $X^2_{8,4}$ follows from that of $X^2_{8,5}$ by considering the exact sequence 
$$\cdots\rightarrow\widetilde{KO}^{-2}(X^2_{8,5})\rightarrow\widetilde{KO}^{-2}(X_{8,4})\rightarrow\widetilde{KO}^{-2}(S^5)\rightarrow\cdots$$
and noting that the last group in the above sequence is zero. Similar considerations shows that $X^2_{8,3}$, $X^2_{8,2}$ are $W$-trivial. This completes the proof in the case $t=1$. 

Next assume that $t>1$ and consider the exact sequence 
$$\cdots\rightarrow\widetilde{KO}^{-2}(X_{8t,8t-3})\stackrel{\alpha} \rightarrow \widetilde{KO}^{-2}(X_{8t,2})\rightarrow \widetilde{KO}^{-2}(X_{8t-3,2})\rightarrow\cdots.$$
The last group in the above sequence is zero (\cite{fujii}, Table 4). Hence $\alpha$ is an epimorphism. We claim that $X^2_{8t,8t-3}$ is $W$-trivial. 
The inclusion map $S^{8t}\longrightarrow X^2_{8t,8t-3}$ induces isomorphism in cohomology in degree $8t$. If $\xi$ is a vector bundle over $X^2_{8t,8t-3}$ with $w(\xi)\neq 1$, then $w_{8t}(\xi)\neq 0$. This bundle then pulls back to a bundle $\eta$ over $S^{8t}$ with $w(\eta)\neq 1$. This is a contradiction as $t>1$. Hence 
$X^2_{8t,8t-3}$ is $W$-trivial. The surjectivity of $\alpha$ now imples that $X^2_{8t,2}$ is $W$-trivial. This completes the proof of (5) and the theorem.

{\em Proof of Theorem\,\ref{fourththeorem}}. First note that if $m=2, 3$, then $X^4_{m,n}$ is always $W$-trivial. Since the obvious map $X^4_{4,n}\longrightarrow S^8$ induces isomorphism in cohomology in top dimension, we have that $X^4_{4,n}$ is not $W$-trivial. This completes the proof of (1). 

Now assume  that $m>4$. Consider the exact sequence 
$$\cdots\rightarrow\widetilde{KO}^{-4}(X_{m,3})\stackrel{\alpha}\longrightarrow \widetilde{KO}^{-4}(\mathbb R\mathbb P^m)\longrightarrow \widetilde{KO}^{-4}(\mathbb R\mathbb P^3)\rightarrow\cdots.$$
The last group in the above sequence is well known to be zero (see \cite{fujii1}) and hence $\alpha$ is an epimorphism. By Theorem\,\ref{tanakastheorem}, $\Sigma^4\mathbb R\mathbb P^m$ is not $W$-trivial. The surjectivity of $\alpha$ now implies that $X^4_{m,3}$ is not $W$-trivial if $m>4$. In view of Proposition\,\ref{firstlemma} (3), 
$X^4_{m,n}$ is not $W$-trivial for $m>4$ and $n=1,2,3$. To complete the proof of the theorem we now show that if $m>4$, then $X^4_{m,4}$ is $W$-trivial. 
To see this consider the exact sequence (see \cite{fujii} Section 3)
$$0\rightarrow \widetilde{KO}^{-4}(X_{m,4})\stackrel{\alpha}\longrightarrow\widetilde{KO}^{-4}(\mathbb R\mathbb P^m)\stackrel{\beta}\longrightarrow\widetilde{KO}^{-4}(\mathbb R\mathbb P^4)\rightarrow 0.$$
The last two groups in the above exact sequence are finite cyclic \cite{fujii1}. By Theorem\,\ref{tanakastheorem}, the spaces $\Sigma^4\mathbb R\mathbb P^m$ and $\Sigma^4\mathbb R\mathbb P^4$ are not $W$-trivial. Thus if $\theta$, $\eta$ are generators of the second and the third group respectively, then we must have $w(\xi)\neq 1$ and $w(\eta)\neq 1$. Assume that $\beta(\theta)=\eta$.  
Now let $\xi \in \widetilde{KO}^{-4}(X_{m,4})$ be such that $w(\xi)\neq 1$. Then $w(\alpha(\xi))\neq 1$ as $\alpha$ induces isomorphism in cohomology in degress $j$, $j\geq 9$. Let $\alpha(\xi)=s\theta$. Then $s$ is odd as the cup products in $\widetilde{H}^*(\Sigma^4\mathbb R\mathbb P^m)$ are all zero. The calculation  
$$w(\beta\alpha(\xi))=w(s\beta(\theta))=w(s\eta)=w(\eta)\neq 1$$
contradicts the exactness of the above sequence. Thus $w(\xi)=1$ for every $\xi \in \widetilde{KO}^{-4}(X_{m,4})$ proving that $X^4_{m,4}$ is $W$-trivial. This shows that $X^4_{m,n}$ is $W$-trivial for all $n\geq 4$ if $m>4$. This completes the proof of the theorem. 

{\em Proof of Theorem\,\ref{fifththeorem}}. Consider the exact sequence 
$$\widetilde{KO}(X_{m,n})\stackrel{\alpha}\longrightarrow\widetilde{KO}(\mathbb R\mathbb P^m)\longrightarrow\widetilde{KO}(\mathbb R\mathbb P^n)\longrightarrow 0.$$
It is known (see \cite{adams}) that the image of $\alpha$ is generated by $2^{\varphi(n)}\xi$ where $\xi$ is the canonical line bundle over $\mathbb R\mathbb P^m$ and $\varphi(n)$ is as in Theorem\,\ref{tanakastheorem1}. It follows that in the exact sequence 
$$\widetilde{KO}^{-8}(X_{m,n})\stackrel{\beta}\longrightarrow\widetilde{KO}^{-8}(\mathbb R\mathbb P^m)\longrightarrow\widetilde{KO}^{-8}(\mathbb R\mathbb P^n)\longrightarrow 0$$
the image of $\beta$ is generated by $2^{\varphi(n)}\eta$ where $\eta$ corresponds to $\xi$ under the Bott periodicity isomorphism. Since $2^{\varphi(n)}$ is even and the cup products in $\widetilde{H}^*(\Sigma^8\mathbb R\mathbb P^m)$ are zero, it follows that $w(2^{\varphi(n)}\eta)=1$. Hence if $\theta$ is in the image of $\beta$, then $w(\theta)=1$. As $\alpha$ induces isomorphism in cohomology in degrees $j$, $j\geq n+9$, it follows that $X^8_{m,n}$ is $W$-trivial. 
This completes the proof of the theorem. 

{\em Acknowledgement}. The second author acknowledges the hospitality of Indian Statistical Institute, Bangalore, where part of this work was done.

\end{document}